\documentclass[12pt]{elsarticle}

\usepackage{amsthm,amsmath,amssymb,
mathrsfs,easymat}
\usepackage[all]{xy}

\usepackage[active]{srcltx}

\DeclareMathOperator{\Ker}{Ker}
\DeclareMathOperator{\tr}{tr}

\renewcommand{\le}{\leqslant}
\renewcommand{\ge}{\geqslant}

\newcommand{\ff}{\mathbb F}
\newcommand{\mc}{\mathcal }

\newcommand{\un}{\underline }

\newcommand{\ind}[1]{\genfrac{}{}{0pt}{}{\scriptstyle \it #1}{}}

\newcommand{\mat}[1]{\begin{bmatrix}
   #1\end{bmatrix}}

\newcommand{\ma}[1]{\begin{matrix}
   #1\end{matrix}}

\newcommand{\arr}[1]{
\left[  \begin{array}#1
\end{array}  \right]}

\newcommand{\wt}{\widetilde}

\newcommand{\s}{\scriptscriptstyle}

\newtheorem{theorem}{Theorem}

\theoremstyle{definition}

\newtheorem{example}{Example}

\theoremstyle{remark}
\newtheorem{remark}{Remark}

\begin{document}
\title{Pairs of commuting nilpotent operators with one-dimensional intersection of kernels and matrices commuting with a Weyr matrix}

\author[ser]{Vitalij M.  Bondarenko} \ead{vit-bond@imath.kiev.ua}

\author[fu]{Vyacheslav Futorny}
\address[fu]{Department of Mathematics, University of S\~ ao Paulo, Brazil}
\ead{futorny@ime.usp.br}

\author[pet]{Anatolii P.~Petravchuk}
\address[pet]{Faculty of Mechanics and Mathematics, Taras Shevchenko University, Kiev, Ukraine}
\ead{aptr@univ.kiev.ua}

\author[ser]{Vladimir~V.~Sergeichuk
}
\ead{sergeich@imath.kiev.ua}
\address[ser]{Institute of Mathematics, Tereshchenkivska 3,
Kiev, Ukraine}

\cortext[cor]{Linear Algebra Appl. 2021,  DOI: 10.1016/j.laa.2020.10.040}

\begin{abstract}
I.M. Gelfand and V.A. Ponomarev (1969) proved
that the problem of classifying pairs $(\mc A,\mc B)$ of commuting nilpotent operators on a vector space contains the problem of classifying an arbitrary $t$-tuple of linear operators. Moreover, it contains the problem of classifying representations of an arbitrary quiver and an arbitrary finite-dimensional algebra, and so it is considered as hopeless.

If $(\mc A,\mc B)$ is such a pair, then $\Ker\mc A\cap\Ker\mc B\ne 0$. We give a simple normal form $(A_{\text{nor}},B_{\text{nor}})$ of the matrices of $(\mc A,\mc B)$ if $\Ker\mc A\cap\Ker\mc B$ is one-dimensional. We do not know whether it is canonical; i.e., whether $(A_{\text{nor}},B_{\text{nor}})$  is uniquely determined by $(\mc A,\mc B)$. We prove its uniqueness only if the Jordan canonical form of $\mc A$ is a direct sum of Jordan blocks of the same size and the field is of zero characteristic.

The matrix $A_{\text{nor}}$ is the Weyr canonical form of $\mc A$, and $B_{\text{nor}}$ commutes with $A_{\text{nor}}$. In order to describe the structure of $(A_{\text{nor}},B_{\text{nor}})$, we describe explicitly all matrices commuting with a given Weyr matrix.


\end{abstract}

\maketitle

\section{Introduction}

We give a simple normal form under similarity of a pair of commuting nilpotent matrices whose common null space is one-dimensional. We also describe explicitly all matrices commuting with a given Weyr matrix.

Gelfand and Ponomarev \cite{gel-pon} proved
that the problem of classifying pairs $(M,N)$ of commuting nilpotent matrices under similarity transformations
\[
(M,N)\mapsto S^{-1}(M,N)S:=(S^{-1}MS,S^{-1}NS),\qquad S\text{ is nonsingular},
\]
contains the problem of classifying $t$-tuples of matrices with any $t$ under similarity transformations
\[
(A_1,\dots,A_t)\mapsto (S^{-1}A_1S,\dots,S^{-1}A_tS),\qquad S\text{ is nonsingular}
\]
(we say that a problem \emph{contains} another problem if solving the first would solve the second).

Moreover, the problem of classifying matrix pairs under similarity contains the problem of classifying each system consisting of vector spaces and linear mappings between them; that is, representations of an arbitrary quiver (see Belitskii and Sergeichuk \cite{bel-ser_compl},  Barot \cite[Section 2.4]{bar1}, and Krause \cite[Section 10]{kra}). It also contains the problem of classifying representations of an arbitrary finite-dimensional algebra
(Barot \cite[Proposition 9.14]{bar}).
Classification problems that contain the problem of classifying matrix pairs under similarity are called \emph{wild}; they are considered as hopeless.

However, two classification results about matrix pairs under similarity were obtained in 1983:

\begin{itemize}
  \item Friedland \cite{fri} gave a system of invariants of matrix pairs with respect to
similarity.
  \item Belitskii \cite{bel} (see also \cite{bel1,ser_can}) constructed an algorithm that reduces by similarity transformations each pair $(M,N)$ of square matrices of the same size over an algebraically closed field to some pair $(M_{\text{can}},N_{\text{can}})$
such that $(M,N)$ is similar to $(P,Q)$ if and only if $(M_{\text{can}},N_{\text{can}})=
(P_{\text{can}},Q_{\text{can}})$.
Thus, $(M_{\text{can}},N_{\text{can}})$
can be considered as a \emph{canonical form} of $(M,N)$ under similarity.
\end{itemize}

In the article \cite{ser_can}, in which the tame-wild theorem is proved using Belitskii's algorithm, Sergeichuk  calls $W_{M}:=M_{\text{can}}$ the \emph{Weyr canonical form} of $M$ since $W_{M}$ is constructed by the Weyr characteristic of $M$. The matrix $W_{M}$  is permutation similar to the Jordan canonical form of $M$ and has the property: all matrices commuting with $W_{M}$ are upper block triangular.
The Weyr canonical matrices are studied in \cite[Section 3.4]{h-j} and \cite{mear}.

Each pair $(M,N)$ is reduced by Belitskii's algorithm \cite{bel}  as follows:
\begin{itemize}
  \item[-] First $(M,N)$ is reduced by similarity transformations to a pair $(W_{M},P)$.
  \item[-] Then $P$ is reduced to its canonical form by those similarity transformations that preserve the Weyr matrix $W_{M}$. For this purpose, Belitskii partitions $P$ conformally with $W_{M}$ and consistently reduces its blocks.
      On each step, Belitskii makes additional partitions into blocks and restricts the set of admissible transformations to those that preserve the already reduced blocks.
\end{itemize}

Let $M$ and $N$ be $n\times n$ commuting nilpotent matrices over a field $\ff$.  Since they commute, $(M,N)$ is similar to a pair $(M',N')$ of upper triangular matrices. Since $M'$ and $N'$ are nilpotent, their main diagonals are zero. Hence, $M'e_1=N'e_1=0$ with $e_1:=(1,0,\dots,0)^T$, and so their common null space
\begin{equation}\label{hbb}
\mc N(M,N):=\{v\in\ff^n\,|\,Mv=Nv=0\}
\end{equation}
is nonzero.

From now on, $(M,N)$ denotes a pair of commuting nilpotent matrices over $\ff$ for which $\dim \mc N(M,N)=1$. We reduce $(M,N)$ by similarity transformations to some simple pair $(W_M,B)$. As in Belitskii's algorithm, $W_M$ is the Weyr canonical form of $M$, we partition $P$ conformally with $W_{M}$ and consistently reduce its blocks by transformations that preserve $W_{M}$. However, we do not use all transformations that preserve $W_{M}$ (unlike Belitskii's algorithm), and so we do not know whether $(W_{M},B)$ is uniquely determined by $(M,N)$. We can prove its uniqueness only if the Jordan canonical form of $M$ is a direct sum of Jordan blocks of the same size and the field $\ff$ is of zero characteristic.

Since $N$ commutes with $M$, $B$ commutes with the Weyr matrix $W_{M}$. In order to describe the structure of the matrix $B$, we describe the form of all matrices that commute with a Weyr matrix.

The article is organized as follows. The main theorem is formulated in
Section \ref{ss1}; it is proved in Section \ref{ss4}.
In Section \ref{ssr} we give a method for constructing all matrices commuting with a given Weyr matrix.

\section{The main theorem}
\label{ss1}

The Weyr canonical form of a nilpotent matrix $M$ over a field $\ff$ is constructed as follows. Let
\begin{equation}\label{dem0}
J=\underbrace{J_{k_1}(0)\oplus\dots\oplus J_{k_1}(0)}_{p_1\text{ summands}}
\oplus\dots\oplus
\underbrace{J_{k_t}(0)\oplus\dots\oplus J_{k_t}(0)}_{p_t\text{ summands}}
\end{equation}
be the Jordan canonical form of $M$,
in which $k_1>k_2>\dots>k_t>0$, $p_1,\dots,p_t$ are nonzero, and
\[
J_k(0):=\mat{0&1&&0\\&0&\ddots
\\&&\ddots&1\\0&&&0}\qquad(\text{$k$-by-$k$}).
\]
Permute the rows of $J$  collecting the first rows of all Jordan blocks at the top, collect the second rows under them, and so on. Make the same permutations of columns of $J$ and obtain a matrix of the form
\begin{equation}\label{als0}
W=\mat{0_{r_1}
&\left[\begin{smallmatrix}I\\0
\end{smallmatrix}\right]
&&0\\
&0_{r_2}
&\ddots&
\\&&\ddots&\left[\begin{smallmatrix}I\\0
\end{smallmatrix}\right]
\\0&&&
0_{r_{k}}
},\qquad k:=k_1.
\end{equation}
The matrix $W$ is called the \emph{Weyr canonical form} of $M$; it is permutation similar to $J$.
The sequence $r_1,\dots,r_{k_1}$ is the \emph{Weyr  characteristics} of $M$; each $r_i$ is the number of Jordan blocks in \eqref{dem0} of sizes equal to or greater than $i\times i$.

Other versions of the following theorem are
given in \cite{bel}, \cite{bel1}, \cite[Theorem 3.4.2.10(a)]{h-j}, \cite[Section 3.2]{mear}, and  \cite[Theorem 1.2]{ser_can}.

\begin{theorem}[Belitskii \cite{bel}]\label{nbq}
The set of matrices commuting with the Weyr matrix \eqref{als0} consists of all the matrices of the form
\begin{equation*}\label{dde}
S=\mat{\un S_{11}&\dots&\un S_{1k}
\\&\ddots&\vdots\\0&&\un S_{kk}},\qquad\text{each $\un S_{ii}$ is }r_i\times r_i,
\end{equation*}
in which
\begin{itemize}
 \item each $\un S_{1j}$ has a staircase form
 \begin{equation}\label{rvt}
\left[\begin{MAT}(c){ccccccc}
&&&&&&*
\\\phantom{A} &&&&&&
\\ \text{\raisebox{-6pt}{$0$}}&&&&\ddots&&\\ &&&&&&
\addpath{(0,3,3)rr}
\addpath{(2,2,3)rr}
\addpath{(5,1,3)rr}\\
\end{MAT}\,\right];
 \end{equation}
all entries on the stairs and over them are arbitrary and all entries under the stairs are zero.
For each $i$, the location of stairs is uniquely determined by $W$, the first stair contains an entry from the first column and the last stair contains the last entry from the last column.

  \item every $\un S_{ij}$ with $1\le i\le j$ is the $r_i\times r_i$ submatrix of $\un S_{1,j-i+1}$ located in its upper left corner.
\end{itemize}
\end{theorem}

A stronger form of Theorem \ref{nbq} is given in Theorem \ref{lww0}.

Our main result is the following theorem, which is proved in Section \ref{ss4}.

\begin{theorem}
\label{t1q0}

\begin{itemize}
  \item[\rm(a)]
Let $(M,N)$ be a pair of commuting nilpotent matrices over a field $\mathbb F$, whose common null space \eqref{hbb} is one-dimensional. Let \eqref{als0} be the Weyr canonical form of $M$.
Then $(M,N)$ is similar to a pair of the form
\begin{equation}\label{bks0}
\setlength{\arraycolsep}{2pt}
(W,B):=\left(\mat{0_{r_1}
&\left[\begin{smallmatrix}I\\0
\end{smallmatrix}\right]
&&0\\
&0_{r_2}
&\ddots&
\\&&\ddots&\left[\begin{smallmatrix}I\\0
\end{smallmatrix}\right]
\\0&&&
0_{r_{k}}
},
\mat{J_{{r_1}}(0)
&\un B_{12}
&\dots&
\un B_{1k}\\
&J_{{r_2}}(0)
&\ddots&\vdots
\\&&\ddots&\un B_{k-1,k}
\\0&&&
J_{r_k}(0)
}
\right)
\end{equation}
such that
\begin{itemize}
  \item its first  matrix $W$ is the Weyr canonical form of $M$,

  \item each $\un B_{1j}$ has the staircase form
 \[
\left[\begin{MAT}(c){ccccccc}
&&&&&&0
\\\phantom{A} &&&&&&
\\ \text{\raisebox{-6pt}{$0$}}&&&&\ddots&&\\ &&&&&&
\addpath{(0,3,3)rr}
\addpath{(2,2,3)rr}
\addpath{(5,1,3)rr}\\
\end{MAT}\,\right]\]
whose stairs are located as the stairs of $\un S_{1j}$ in \eqref{rvt}, and all nonzero entries of $\un B_{1j}$ are located on the stairs,

  \item each $\un B_{1+l,j+l}$ with $l\ge 1$
is a submatrix of
$\un B_{1j}$
that is located in its upper left corner.

\end{itemize}

  \item[\rm(b)] In particular, if the Jordan canonical form of $M$ is
      \begin{equation*}\label{54d}
J:=J_{k}(0)\oplus\dots\oplus J_{k}(0)\qquad (r\text{ summands}),
\end{equation*}
then $(M,N)$ is similar to a pair of the form
\begin{equation}\label{b6s0}
\setlength{\arraycolsep}{2pt}
(W,B):=\left(\mat{0_r&I_r&&&0_r\\
&0_r&I_r\\
&&0_r&\ddots\\&&&\ddots&I_r\\0_r&&&&0_r},
\mat{J_r(0)&B_1&B_2&\ddots&B_{k-1}\\
&J_r(0)&B_1&\ddots&\ddots\\
&&J_r(0)&\ddots&B_2\\
&&&\ddots&B_1\\
0_r&&&&J_r(0)
}
\right)
\end{equation}
in which every $B_i$ is a matrix whose last row is arbitrary and the other rows are zero.
If the characteristic of\/ $\mathbb F$ is zero, then the pair \eqref{b6s0} is uniquely determined by $(M,N)$.
\end{itemize}
\end{theorem}

\begin{remark}\label{vdi}
The fact that the pair \eqref{b6s0} over a field of zero characteristic  is uniquely determined by $(M,N)$ is unexpected since the similarity transformations with $(W,B)$ that preserve all entries except for the entries denoted by stars have free parameters. Our proof of (b) is technical and does not explain why these parameters do not change the entries denoted by stars.
We do not know whether the pair \eqref{bks0} is uniquely determined by $(M,N)$ if the Jordan canonical form of $M$ has Jordan blocks of unequal sizes.
\end{remark}

\begin{example}\label{1n3x}
Let $(M,N)$ be a pair of commuting nilpotent matrices whose common null space is one-dimensional. Let the Jordan canonical form of $M$ be one of the matrices
\begin{gather}\label{11a}
\underbrace{\oplus J_{4}(0)}_{p\text{ summands}}
\oplus \underbrace{\oplus J_{3}(0)}_{q\text{ summands}}
\oplus \underbrace{\oplus J_{2}(0)}_{r\text{ summands}}
\oplus
\underbrace{\oplus J_{1}(0)}_{s\text{ summands}},\\ \label{11b}
\underbrace{J_7(0)\oplus\dots\oplus J_7(0)}_{p\text{ summands}}\oplus
\underbrace{J_4(0)\oplus\dots\oplus J_4(0)}_{q\text{ summands}}\oplus
\underbrace{J_2(0)\oplus\dots\oplus J_2(0)}_{r\text{ summands}},
\end{gather}
in which $p,q,r,s$ are nonzero.
In Example \ref{1n3}, we describe
all matrices that commute with the Weyr canonical forms of \eqref{11a} and \eqref{11b}.
Using it and Theorem \ref{t1q0}(a), we obtain the following normal forms for $(M,N)$.

(a) If the Jordan matrix of $M$ is \eqref{11a}, then $(M,N)$ is similar to
\[
\setlength{\arraycolsep}{2pt}
\left(\left[\begin{MAT}(@){c1c1c1c}
\ma{0_{p}&&&\\
&0_{q}\\ &&0_{r}\\
&&&0_{s}\\}
                   &
\ma{I_{p}&0&0\\
0&I_{q}&0\\
0&0&I_{r}\\
0&0&0\\}
&0&0
                            \\1
0&\ma{0_{p}&&\\
&0_{q}\\ &&0_{r}\\}
&
\ma{I_{p}&0\\
0&I_{q}\\
0&0\\}
&0
        \\1
0&0&\ma{0_{p}\\
&0_{q}\\}&\ma{I_{p}\\
0\\}
      \\1
0&0&0&0_{p}
\\
\end{MAT}\right] ,\
\left[\begin{MAT}(@){c1c1c1c}
 J_{p+q+r+s}\
&\ma{0&0&0\\
B_2&0&0\\ 0&B_3&0\\
0&0&B_4\\}&\ma{0&0\\
0&0\\ C_3&0\\
0&C_4\\}&\ma{0\\
0\\ 0\\
D_4\\}
             \\1
0&J_{p+q+r}&\ma{0&0\\
B_2&0\\ 0&B_3\\}&\ma{0\\
0\\ C_3\\}
        \\1
0&0&J_{p+q}&\ma{0\\
B_2\\}
      \\1
0&0&0&J_{p}
\\
\end{MAT}\right]\right)
\]
in which $J_m:=J_m(0)$ for all $m$
and the blocks $B_2, B_3, B_4, C_3, C_4, D_4$
 have the form
\begin{equation}\label{1a20}
 \mat{0&0&\dots&0\\[-5pt]\hdotsfor{4} \\0&0&\dots&0\\ *& *&\dots&*},
\end{equation}
where the stars denote arbitrary elements of\/ $\mathbb F$.

(b) If the Jordan matrix of $M$ is \eqref{11b}, then $(M,N)$ is similar to a pair of the form
\eqref{bks0}, in which
 \begin{equation*}\label{ebr1}
[\un B_{12}|\dots|\un B_{17}]=\arr{{ccc|cc|cc|c|c|c}
\Delta&0&0 &\Delta&0 & 0&0 & 0 & 0 & 0\\
0&\Delta&0 & 0&0 &\Delta&0 &\Delta & 0 & 0\\
0&0&\Delta & 0&\Delta & 0&\Delta & 0 & \Delta & \Delta\\
},
\end{equation*}
the sizes of horizontal strips are $p,q,r$, the sizes of vertical strips are $p,q,r,p,q,p,q,p,p,p$, and  each $\Delta$ denotes a block of the form \eqref{1a20}.
\end{example}

\section{Matrices commuting with a Weyr matrix}\label{ssr}

The most important property of each Weyr canonical matrix $W$ was found by Belitskii: all matrices $S$ commuting with $W$ are block triangular. However, all versions of
Belitskii's theorem in \cite{bel,bel1,h-j,mear,ser_can} do not give the positions of stairs in \eqref{rvt}. We give them in this section, which is important for constructing of pairs \eqref{bks0} (see Example \ref{1n3x}) and for other applications of Belitskii's theorem.

\subsection{Belitskii's theorem}

A \emph{Weyr matrix} over a field $\ff$ is each matrix of the form
\begin{equation}\label{5dr}
W=(\lambda _1I_{n_1}+W_1)\oplus\dots\oplus (\lambda _lI_{n_l}+W_l),
\end{equation}
in which $\lambda _1,\dots,\lambda _l\in \ff$ are distinct and $W_1,\dots,W_l$
are nilpotent Weyr matrices of the form \eqref{als0}.
Each matrix over $\ff$ in which all with distinct eigenvalues $\lambda _1,\dots,\lambda _l$ are contained in $\ff$ is similar to some matrix \eqref{5dr}. If a matrix $S$ commutes with $W$, then it has the form
\[
S=S_1\oplus\dots\oplus S_l,
\]
in which $S_1,\dots,S_l$ commute with $W_1,\dots,W_l$, respectively. Therefore, it suffices to describe all matrices commuting with a nilpotent Weyr matrix.

Let $J$ be the Jordan matrix \eqref{dem0}. Let us consider the set of its Jordan blocks of sizes  equal to or greater than $i\times i$, and let  $m_i$ be the number of Jordan blocks of distinct sizes in this set; that is,
\begin{equation}\label{dds0}
(m_1,\dots,m_{k})=\bigr(\underbrace{t,\dots,
t}_{k_t\text{ times}},\
\underbrace{t-1,\dots,
t-1}_{k_{t-1}-k_t\text{ times}},\ \dots,\!\!
\underbrace{1,\dots,
1}_{k-k_2\text{ times}}\bigl).
\end{equation}

\begin{theorem}\label{lww0}
Let $W$ in \eqref{als0} be the Weyr canonical form of the Jordan matrix \eqref{dem0}.
Then the set of matrices that commute with $W$ consists of all matrices of the form
\begin{equation}\label{xkx0}
S=\left[\begin{MAT}(r){c1c1c1c}
\un S_{11}&\un S_{12}&\ddots &\un S_{1k}\\1
&\un S_{22}&\ \ddots\ &\ddots\\1
&&\ddots&\un S_{\,k-1,k}\\1
0&&&\un S_{\,kk}
\\\end{MAT}\right],\qquad
\un S_{\,\alpha \beta }=
\mat{S_{11}^{(\beta -\alpha)}&\dots&S_{1m_{\beta }}^{(\beta -\alpha) }
\\ \hdotsfor{3}\\
S_{m_{\alpha }1}^{(\beta -\alpha)}&\dots&S_{m_{\alpha }m_{\beta }}^{(\beta -\alpha) }}
\end{equation}
($m_1,m_2,\dots$ are defined in \eqref{dds0}),
in which every $S_{ij}^{(\gamma )}$ is a $p_i\times p_j$ matrix,
\begin{equation}\label{xsw0}
\begin{array}{ll}
S_{ij}^{(\gamma) }\text{ is arbitrary}&\text{ \ if \ }
k_j-k_i\le\gamma,
     \\
S_{ij}^{(\gamma) }=0&\text{ \ if \ }
k_j-k_i>\gamma.
\end{array}
\end{equation}
\end{theorem}

\begin{proof}
The nilpotent Jordan matrix $J$ in
\eqref{dem0} is permutation similar to the matrix
\begin{equation*}\label{jdj}
J^{\text{\it +}}:=J_{k_1}(0_{p_1})\oplus
\dots\oplus J_{k_t}(0_{p_t}),\qquad
k_1>k_2>\dots>k_t>0,
\end{equation*}
in which
\[
J_k(0_p):=\mat{0_p&I_p&&0_p\\&0_p&\ddots
\\&&\ddots&I_p\\0_p&&&0_p}\qquad(\text{$k^2$ blocks}).
\]

By Gantmacher \cite[Chapter VIII, \S\,2]{gan},
each matrix that commutes with $J^{\text{\it +}}$ has the form
\[
\newcommand{\vvdots}%
{\text{\raisebox{3pt}{$\vdots$}}}
\newcommand{\vddots}%
{\text{\raisebox{3pt}{$\ddots$}}}
C^{\text{\it +}}\!\!=\begin{MAT}(@){cccccccccccccl}
S_{11}^{\s (0)}&S_{11}^{\s (1)}&S_{11}^{\s (2)}
&\ddots&S_{11}^{\s (k_1-1)}
      &
S_{12}^{\s(0)}&S_{12}^{\s(1)}
&\ddots&S_{12}^{\s(k_2-1)}
     &\dots&
S_{1t}^{\s(0)}&\ddots&S_{1t}^{\s(k_t-1)}
                       &  \\  
&S_{11}^{\s(0)}&S_{11}^{\s(1)}
&\ddots&S_{11}^{\s(k_1-2)}
      &
&S_{12}^{\s(0)}&\ddots&\ddots
     &\dots&
&\ddots&\ddots
                       &  \\  
&&S_{11}^{\s(0)}&\ddots&\ddots
      &
&&\ddots&S_{12}^{\s(1)}
     &\dots&
&&S_{1t}^{\s(0)}
             & \scriptstyle k_1 \\  
&&&\ddots&S_{11}^{\s(1)}
      &
&&&S_{12}^{\s(0)}
     &\dots&
&&
                       &  \\  
&&&&S_{11}^{\s(0)}
      &
&&&
     &\dots&
&&
                      &     \\
&\!\!\!\!S_{21}^{\s(k_1-k_2)}&\!\!S_{21}^{\s(k_1-k_2+1)}\!\!
&\ddots\!\!&S_{21}^{\s(k_1-1)}
      &
S_{22}^{\s(0)}&S_{22}^{\s(1)}
&\ddots&S_{22}^{\s(k_2-1)}
     &\dots&
S_{2t}^{\s(0)}&\ddots&S_{2t}^{\s(k_t-1)}
                       & \\ 
&&S_{21}^{\s(k_1-k_2)}
&\ddots&S_{21}^{\s(k_1-2)}
      &
&S_{22}^{\s(0)}&\ddots&\ddots
         &\dots&
&\ddots&\ddots
     & \scriptstyle k_2 \\ 
&&&\ddots&\ddots
      &
&&\ddots&S_{22}^{\s(1)}
     &\dots&
&&S_{2t}^{\s(0)}
                      & \\ 
&&&&S_{21}^{\s(k_1-k_2)}
      &&&
&S_{22}^{\s(0)}
     &\dots&
&&
                        &   \\ 
\vvdots&\vvdots&\vvdots&\vvdots
&\vvdots
&\vvdots&\vvdots&\vvdots
&\vvdots&\vddots&\vvdots
&\vvdots&\vvdots
                       &    \\
&&S_{t1}^{\s(k_1-k_t)}&\ddots&S_{t1}^{\s(k_1-1)}
      &
&\!\!\!\!\!S_{t2}^{\s(k_2-k_t)}\!\!
&\ddots&S_{t2}^{\s(k_2-1)}
     &\dots&
S_{tt}^{\s(0)}&\ddots&S_{tt}^{\s(k_t-1)}
                      &    \\ 
&&&\ddots&\ddots
      &
&&\ddots&\ddots
     &\dots&
&\ddots&\ddots
  & \scriptstyle k_t    \\ 
&&&&S_{t1}^{\s(k_1-k_t)}
      &&&
&\!\!\!S_{t2}^{\s(k_2-k_t)}
     &\dots&
&&S_{tt}^{\s(0)}
          &    \\
&&\scriptstyle k_1&&&&&\scriptstyle k_2&&&&\scriptstyle k_t\!\!\!\!\!\!\!\!&&
\addpath{(0,1,4)rrrrrrrrrrrrr%
uuuuuuuuuuuuulllllllllllllddddddddddddd}
\addpath{(5,1,1)uuuuuuuuuuuuu}
\addpath{(9,1,1)uuuuuuuuuuuuu}
\addpath{(10,1,1)uuuuuuuuuuuuu}
\addpath{(0,4,1)rrrrrrrrrrrrr}
\addpath{(0,5,1)rrrrrrrrrrrrr}
\addpath{(0,9,1)rrrrrrrrrrrrr}
          \\
          \end{MAT}
\]
The sets of rows of $J^{\text{\it +}}$ and $C^{\text{\it +}}$ are partitioned into $t$ superstrips, each of them is partitioned into $k_1, k_2, \dots, k_t$ strips, respectively. As in \cite[Section 1.3]{ser_can},
we permute the horizontal and vertical strips in $J^{\text{\it +}}$ and $C^{\text{\it +}}$, gathering the first strips of each superstrip, then the second strips of each superstrip, and so on, until we obtain  $W$ of the form \eqref{als0} and $S$ of the form \eqref{xkx0}.

It is observed from the picture that the blocks of $C^{\text{\it +}}$ satisfy \eqref{xsw0}.
Hence, the blocks of $S$ satisfy \eqref{xsw0} too.
\end{proof}

Note that each superblock  $\un S_{\,\alpha \beta }$ in  \eqref{xkx0} consists of the following blocks of $C^{\text{\it +}}$:
\emph{the $(i,j)$th block $S_{ij}^{(\beta-\alpha)}$ of $\un S_{\,\alpha \beta }$ is the $(\alpha, \beta)$th block of the $(i,j)$th superblock of $C^{\text{\it +}}$}.

\begin{example}\label{aak} The nilpotent Jordan matrix
\[
J=\underbrace{J_3(0)\oplus\dots\oplus J_3(0)}_{p\text{ summands}}\oplus
\underbrace{J_2(0)\oplus\dots\oplus J_2(0)}_{q\text{ summands}}\oplus
\underbrace{J_1(0)\oplus\dots\oplus J_1(0)}_{r\text{ summands}}
\]
is permutation similar to $J^{\text{\it +}}=J_3(0_p)\oplus J_2(0_q)\oplus 0_r$.
This matrix and each matrix $C^{\text{\it +}}$ that commutes with it have the form
\[
\ \begin{matrix}
J^{\text{\it +}}= \\ \\[-7pt]
\end{matrix}\
\begin{MAT}(@){ccccccl}
\,0_p\,&I_p&&   &&&\scriptstyle \it 1\\
&0_p&I_p&   &&&\scriptstyle \it 2\\
&&0_p&   &&&\scriptstyle \it 3\\
&&&0_q    &I_q&&\scriptstyle \it 1'\\
&&&   &0_q&&\scriptstyle \it 2'\\
&&&   &&0_r&\scriptstyle \it 1''\\
\scriptstyle \it 1&\scriptstyle
\it 2&
\scriptstyle \it 3&\scriptstyle
\it 1'&
\scriptstyle \it 2'&\scriptstyle \it 1''&
\addpath{(0,1,4)rrrrrruuuuuulllllldddddd}
\addpath{(3,1,1)uuuuuu}
\addpath{(5,1,1)uuuuuu}
\addpath{(0,2,1)rrrrrr}
\addpath{(0,4,1)rrrrrr}\\
\end{MAT}\quad
\begin{matrix}
C^{\text{\it +}}=\,\\ \\[-5pt]
\end{matrix}\
\begin{MAT}(@){ccccccl}
S_{11}^{(0)}&S_{11}^{(1)}&S_{11}^{(2)}&S_{12}^{(0)}   &S_{12}^{(1)}&S_{13}^{(0)}&\scriptstyle \it 1\\
&S_{11}^{(0)}&S_{11}^{(1)}&  &S_{12}^{(0)}&
&\scriptstyle \it 2\\
&&S_{11}^{(0)}& &&&\scriptstyle \it 3\\
&S_{21}^{(1)}&S_{21}^{(2)}&S_{22}^{(0)}    &S_{22}^{(1)}&S_{23}^{(0)}&\scriptstyle \it 1'\\
&&S_{21}^{(1)}&   &S_{22}^{(0)}&&\scriptstyle \it 2'\\
&&S_{31}^{(2)}&  &S_{32}^{(1)}&S_{33}^{(0)}&\scriptstyle \it 1''\\
\scriptstyle \it 1&\scriptstyle
\it 2&
\scriptstyle \it 3&\scriptstyle
\it 1'&
\scriptstyle \it 2'&\scriptstyle
\it 1''&
\addpath{(0,1,4)rrrrrruuuuuulllllldddddd}
\addpath{(3,1,1)uuuuuu}
\addpath{(5,1,1)uuuuuu}
\addpath{(0,2,1)rrrrrr}
\addpath{(0,4,1)rrrrrr}\\
\end{MAT}
\]
The set of rows of $J^{\text{\it +}}$ is partitioned into 3 \emph{superstrips}, each of them is partitioned into 3, 2, and 1 strips, respectively. We permute the horizontal strips of $J^{\text{\it +}}$ and of $C^{\text{\it +}}$ gathering atop the first strips of each superstrip, then the second strips, and finally the third strips. Making the same permutations of vertical strips, we obtain
\[
\ \begin{matrix}
J^{\#}= \\ \\ \\
\end{matrix}\
\begin{MAT}(@){ccccccl}
\,0_p\,&&&I_p   &&&\scriptstyle \it 1\\
&0_q&&   &I_q&&\scriptstyle \it 1'\\
&&0_r&   &&&\scriptstyle \it 1''\\
&&&0_p    &&I_p&\scriptstyle \it 2\\
&&&   &0_q&&\scriptstyle \it 2'\\
&&&   &&0_p&\scriptstyle \it 3\\
\ind{1}&\ind{1'}&
\ind{1''}&\ind{2}&
\ind{2'}&\ind{3}&
\addpath{(0,1,4)rrrrrruuuuuulllllldddddd}
\addpath{(3,1,1)uuuuuu}
\addpath{(5,1,1)uuuuuu}
\addpath{(0,2,1)rrrrrr}
\addpath{(0,4,1)rrrrrr}\\
\end{MAT}\quad
\begin{matrix}
C^{\#}=\,\\ \\[10pt]
\end{matrix}\
\begin{MAT}(@){ccccccl}
S_{11}^{(0)}&S_{12}^{(0)}&S_{13}^{(0)}&S_{11}^{(1)}   &S_{12}^{(1)}&S_{11}^{(2)}&\scriptstyle \it 1\\
&S_{22}^{(0)}&S_{23}^{(0)}&S_{21}^{(1)}  &S_{22}^{(1)}&S_{21}^{(2)}
&\scriptstyle \it 1'\\
&&S_{33}^{(0)}& &S_{32}^{(1)}&S_{31}^{(2)}&\scriptstyle \it 1''\\
&&&S_{11}^{(0)}    &S_{12}^{(0)}&S_{11}^{(1)}&\scriptstyle \it 2\\
&&&   &S_{22}^{(0)}&S_{21}^{(1)}&\scriptstyle \it 2'\\
&&&  &&S_{11}^{(0)}&\scriptstyle \it 3\\
\ind{1}&\ind{1'}&
\ind{1''}&\ind{2}&
\ind{2'}&\ind{3}&
\addpath{(0,1,4)rrrrrruuuuuulllllldddddd}
\addpath{(3,1,1)uuuuuu}
\addpath{(5,1,1)uuuuuu}
\addpath{(0,2,1)rrrrrr}
\addpath{(0,4,1)rrrrrr}\\
\end{MAT}
\]
in which $W:=J^{\#}$ is the Weyr canonical form of $J$. Each matrix commuting with $J^{\#}$ has the form $C^{\#}$. The matrix $C^{\#}$ satisfies \eqref{xsw0} since $(k_1,k_2,k_3)=(3,2,1)$.
\end{example}

\subsection{Method of constructing all matrices that commute with a given Weyr matrix}
\label{bbtt}

The matrix \eqref{xkx0} can be constructed
by the numbers $k_{\alpha }$ and $m_{\alpha }$ (defined in \eqref{dem0} and \eqref{dds0}) as follows:
\begin{itemize}
  \item First we construct the skew-symmetric matrix
\[
K:=\mat{k_1-k_1&k_2-k_1&\dots&k_t-k_1\\
k_1-k_2&k_2-k_2&\dots&k_t-k_2
\\ \hdotsfor{4}\\
k_1-k_t&k_2-k_t&\dots&k_t-k_t}.
\]

  \item
Then we construct the block matrix
\begin{equation*}\label{ove}
H:=[H_1|H_2|\dots|H_{k_1}],\qquad\text{each $H_{\beta  }$ is }m_1\times m_{\beta},
\end{equation*}
in which the $j$th column of $H_{\beta }$ is obtained from the $j$th column of $K$ by replacing all its entries that are $<\beta $ by
the multiplication sign $\times$ and the other entries by $0$.

  \item Each $\un S_{1 \beta }$ in \eqref{xkx0} is obtained by replacing in $H_{\beta }$ all $\times$'s  by arbitrary blocks of appropriate sizes and all zero entries by zero blocks.

  \item
Each $\un S_{1+l, \beta+l }$ with $l\ge 1$ is an \emph{angular submatrix} of  $\un S_{1 \beta}$; i.e., it is located in the upper left corner of $\un S_{1\beta }$.
\end{itemize}

\begin{example}\label{1n3}
Let us describe all matrices that commute with the Weyr matrices whose Jordan forms are \eqref{11a} and \eqref{11b}.

(a)
Let $J$ be the Jordan matrix \eqref{11a}.
Then $(p_1,p_2,p_3,p_4)=(p,q,r,s)$,
\[
(k_1,k_2,k_3,k_4)=(4,3,2,1),\qquad
(m_1,m_2,m_3,m_4)=(4,3,2,1),
\]
and so
\[
 K=\mat{0&-1&-2&-3
\\1&0&-1&-2\\2&1&0&-1\\3&2&1&0},
         \quad
H=\arr{{cccc|ccc|cc|c}
\times&\times&\times  & \times&\times&\times&\times  & \times&\times & \times\\
0&\times&\times  & \times&
\times&\times&\times  & \times&\times & \times\\
0&0&\times  & \times&
0&\times&\times  & \times&\times & \times\\
0&0&0& \times&0&0&\times & 0&\times & \times\\
}.
\]
Hence the Weyr canonical form of $J$ and an arbitrary matrix commuting with it have the form
\[
\setlength{\arraycolsep}{2pt}
\left[\begin{MAT}(@){c1c1c1c}
\ma{0_{p}&&&\\
&0_{q}\\ &&0_{r}\\
&&&0_{s}\\}
                   &
\ma{I_{p}&0&0\\
0&I_{q}&0\\
0&0&I_{r}\\
0&0&0\\}
&0&0
                            \\1
0&\ma{0_{p}&&\\
&0_{q}\\ &&0_{r}\\}
&
\ma{I_{p}&0\\
0&I_{q}\\
0&0\\}
&0
        \\1
0&0&\ma{0_{p}\\
&0_{q}\\}&\ma{I_{p}\\
0\\}
      \\1
0&0&0&0_{p}
\\
\end{MAT}\right] ,\
\left[\begin{MAT}(@){c1c1c1c}
\ma{A_1&A_1'&A_1''&A_1'''\\
0&A_2&A_2'&A_2''\\
0&0&A_3&A_3'\\
0&0&0&A_4}
&\ma{B_1'&B_1''&B_1''\\
B_2&B_2'&B_2''\\
0&B_3&B_3'\\
0&0&B_4\\}&
\ma{C_1''&C_1'''\\
C_2'&C_2''\\ C_3&C_3'\\
0&C_4\\}&
\ma{D_1'''\\
D_2''\\ D_3'\\
D_4\\}
             \\1
0&
\ma{A_1&A_1'&A_1''\\
0&A_2&A_2'\\
0&0&A_3}
&
\ma{B_1'&B_1''\\
B_2&B_2'\\
0&B_3}
&\ma{C_1''\\
C_2'\\ C_3}
        \\1
0&0&
\ma{A_1&A_1'\\
0&A_2}
&\ma{B_1'\\
B_2\\}
      \\1
0&0&0&A_1
\\
\end{MAT}\right]
\]
in which all blocks $A_i^{(j)},B_i^{(j)}, C_i^{(j)}, D_i^{(j)}$ are arbitrary.

(b) Let $J$ be the Jordan matrix \eqref{11b}.
Then $t=3$,
\[
(k_1,k_2,k_3)=(7,4,2),\qquad
(m_1,\dots,m_7)=(3,3,2,2,
1,1,1),
\]
\[ K=\mat{7-7&4-7&2-7
\\7-4&4-4&2-4\\7-2&4-2&2-2}
=\mat{0&-3&-5\\3&0&-2\\5&2&0},\]
\begin{equation*}\label{ebr}
H=[H_1|\dots|H_7]
=\arr{{ccc|ccc|cc|cc|c|c|c}
\times&\times&\times &
\times&\times&\times & \times&\times & \times&\times &\times &\times &\times\\
0&\times&\times &
0&\times&\times & 0&\times & \times&\times & \times & \times& \times\\
0&0&\times &
0&0&\times & 0&\times & 0&\times & 0 & \times & \times\\
}.
\end{equation*}
\end{example}


\section{Proof of Theorem \ref{t1q0}}\label{ss4}

\subsection{Proof of Theorem \ref{t1q0}(a)}

Each pair $(M,N)$ of commuting nilpotent matrices  whose common null space is one-dimensional is similar to $(W,B)$, in which $W$ is the Weyr canonical form \eqref{als0} of $M$. By Theorem \ref{lww0}, \begin{equation*}\label{wsr}
B=\left[\begin{MAT}(r){c1c1c1c}
\un B_{11}&\un B_{12}&\ \ddots\ &\un B_{1k}\\1
&\un B_{22}&\ddots&\ddots\\1
&&\ddots&\un B_{\,k-1,k}\\1
0&&&\un B_{\,kk}
\\\end{MAT}\right]
                     ,\qquad
\un B_{\alpha \beta }=
\mat{B_{11}^{(\beta-\alpha )}&\dots &B_{1m_{\beta }}^{(\beta-\alpha ) }
\\ \hdotsfor{3}\\
B_{m_{\alpha }1}^{(\beta -\alpha)}&\dots&B_{m_{\alpha }m_{\beta }}^{(\beta -\alpha) }},
\end{equation*}
in which every $B_{ij}^{(\gamma )}$ is an arbitrary $p_i\times p_j$ matrix such that
$B_{ij}^{(\gamma)}=0$ if
$k_j-k_i>\gamma$.

We reduce $(W,B)$ by those similarity transformations $S^{-1}(W,B)S$ that preserve $W$. Then
$WS=SW$, and so $S$ is of the form \eqref{xkx0}.

Since $B$ and $S$ commute with $W$, they are conformally partitioned and have the same staircase form.

Since $B$ is reduced by transformations $S^{-1}BS$, the superblock $\un B_{11}$ is reduced by similarity transformations $\un S_{11}^{-1}\un B_{11}\un S_{11}$, in which
\[
\un B_{11}=\mat{B_{11}^{(0)}&B_{12}^{(0)}
&\ddots&B_{1m_1}^{(0)}\\
&B_{22}^{(0)}&\ddots&\ddots\\
&&\ddots&B_{m_1-1,m_1}^{(0)}\\
0&&&B_{m_1m_1}^{(0)}
},\quad
\un S_{11}=\mat{S_{11}^{(0)}&S_{12}^{(0)}&\ddots
&S_{1m_1}^{(0)}\\
&S_{22}^{(0)}&\ddots&\ddots\\
&&\ddots&S_{m_1-1,m_1}^{(0)}\\
0&&&S_{m_1m_1}^{(0)}
}.
\]
We reduce each $B_{ii}^{(0)}$ to the Jordan canonical form by transformations $\bigl(S_{ii}^{(0)}\bigr)^{-1}B_{ii}^{(0)}
S_{ii}^{(0)}$. The matrix $B$ is nilpotent, hence each $B_{ii}^{(0)}$ is nilpotent too, and so each $B_{ii}^{(0)}$ is a direct sum of singular Jordan blocks.
Thus, $\un B_{11}$ is of the form
\[
\un B_{11}=\mat{0&b_{12}&\ddots&b_{1r_1}\\
&0&\ddots&\ddots\\
&&\ddots&b_{r_1-1,r_1}\\
0&&&0
}.
\]
All $b_{i,i+1}$ are nonzero since the common null space of $W$ and $B$ is one-dimensional.

We reduce $\un B_{11}$ to $J_{r_1}(0)$ by upper triangular similarity transformations (which are admissible since $\un S_{11}$ is upper block-triangular) as follows. First we multiply the second column by $b_{12}^{-1}$ and the second row by $b_{12}$; we obtain $b_{12}=1$. Then we make zero each $b_{1j}$ with $j\ge 3$ by adding the second column multiplied by $-b_{1j}$ to the $j$th column; the inverse transformations of rows change the second row. We obtain $\un B_{11}$ with the first row $(0,1,0,\dots,0)$. In the same way, we make $b_{23}=1$ and $b_{24}=b_{25}=\dots=0$ (the inverse transformations of rows change the third row) and so on, until we obtain $\un B_{11}=J_{r_1}(0)$.

By \eqref{xkx0},  each $\un B_{\,\alpha \alpha }$ is an angular submatrix of $\un B_{11}$. Hence $\un B_{\,\alpha \alpha }=J_{r_\alpha}(0)$ for all $\alpha$.

The $\beta $th  \emph{upper superdiagonal} of $B$ is the sequence of superblocks
\[
(\un B_{1,\beta+1 },\un B_{2,\beta +2},\dots,\un B_{\,k_1-\beta,k_1}),\qquad 0\le \beta \le k_1-1.
\]
Reasoning by induction, we fix $\beta \ge 1$, assume that if $\beta \ge 2$ then the 1st,\ \dots,\ $(\beta -1)$st upper superdiagonals  of $B$ have the form described in Theorem \ref{t1q0}(a), and reduce  the $\beta$th upper superdiagonal to this form as follows. Let $S=I-G$ be of the form \eqref{xkx0}, in which all entries of $G$ are zero except for the entries of $\beta$th upper superdiagonal $(\un G_{1,\beta+1 },\dots,\un G_{\,k_1-\beta,k_1})$. Then the transformation
\begin{align*}
B\mapsto
S^{-1}BS&=(I+G+G^2+\cdots)B(I-G)
\\&=B+(GB-BG)+G(GB-BG)+\cdots
\end{align*}
preserves the 1st,\ \dots,\ $(\beta -1)$st upper superdiagonals of $B$ and reduces the superblock $\un B_{1,\beta+1}$ by the transformation
\begin{equation}\label{1230}
\un B_{1,\beta+1}\mapsto\wt{\un B}_{1,\beta+1}:=
\un B_{1,\beta +1}+\un G_{1,\beta+1 }J_{r_{\beta+1}}(0)-J_{r_{1}}(0)\un G_{1,\beta+1}.
\end{equation}

By Theorem \ref{lww0}, some of the entries of $\un B_{1,\beta+1}$ and the entries of $\un G_{1,\beta+1}$ at the same positions are arbitrary (they are the entries of those blocks that correspond to $\times$'s in $ H_{\beta+1}$; see Section \ref{bbtt}). The other entries of $\un B_{1,\beta+1}$ and $\un G_{1,\beta+1}$ are zero; we denote them by $\emptyset$ (they are the entries of the zero blocks that correspond to 0's in $H_{\beta+1}$).
All entries above the main diagonal are arbitrary. If an entry is arbitrary, then all entries to the right and above it are arbitrary too.
Thus, $\un B_{1,\beta+1}$ and $\un G_{1,\beta+1}$ are staircase matrices with the same location of stairs.

\begin{description}
  \item[Step 1,] \emph{in which we reduce the entries of $\un B_{1,\beta+1}$ on its main diagonal and under it.}

The $\ell$th  \emph{lower diagonal} of $\un B_{1,\beta+1 }=[b_{ij}]_{i=1}^{r_1}{}_{j=1}^{r_{\beta+1 }}$ is the sequence of its entries
\begin{equation*}\label{aasq}
B_{-(\ell)}:=(b_{\ell+1, 1},b_{\ell+2,2},\dots),\qquad
0\le\ell\le r_1-1.
\end{equation*}
Reasoning by induction, we assume that all the lower diagonals of $\un B_{1\beta }$ under $B_{-(\ell)}$ have the form described in Theorem \ref{t1q0}(a). Let $B_{-(\ell)}$ do not have this form. Then $B_{-(\ell)}\ne (\emptyset,\dots,\emptyset)$. Let $G_{-(\ell+1)}:=(g_{\ell+2, 1},g_{\ell+3,2},\dots)$ be the $(\ell+1)$st lower diagonal
 of $\un G_{1,\beta+1 }=[g_{ij}]_{i=1}^{r_1}{}_{j=1}^{r_{\beta+1 }}$. Each
transformation \eqref{1230} changes $B_{-(\ell)}$ as follows: \begin{equation}\label{2340}
\wt{B}_{-(\ell)}=B_{-(\ell)} +(0,g_{\ell+2, 1},g_{\ell+3,2},\dots)-(g_{\ell+2, 1},g_{\ell+3,2},g_{\ell+4,3},\dots).
\end{equation}

Let us
consider a fragment $g_1, g_2, \dots,g_k$ of $G_{-(\ell+1)}$ of the form
\begin{equation}\label{vbv0}
G_{-(\ell+1)}
=(\,\underbrace{\text{any},
\emptyset}_{p\text{ entries}}\, ,\underbrace{g_1, g_2, \dots,g_k}_{k\text{ entries}},  \underbrace{\emptyset,
\text{any}}_{q\text{ entries}}\,),\qquad\text{all }g_i\ne\emptyset,
\end{equation}
in which $p,q\in\{0,1,\dots\}$ and $k\ge 1$ (which means that $g_k$ is on a stair and $g_1, g_2, \dots,g_{k-1}$ are over this stair).

\begin{itemize}
  \item
Let $g_k$ be not in the last column of $\un B_{1,\beta+1}$ without its last entry. By
\eqref{2340},
\begin{equation*}\label{fgh}
\wt B_{-(\ell)}=B_{-(\ell)} +(\underbrace{\dots}_{p\text{ entries}} ,\underbrace{-g_1, g_1-g_2, \dots,g_{k-1}-g_k,g_k}_{k+1\text{ entries}},\underbrace{\dots}_{q\text{ entries}}).
\end{equation*}
We choose $g_1, g_2, \dots,g_k$ in \eqref{vbv0} such that
\begin{equation*}\label{ccd}
\wt B_{-(\ell)}=(\underbrace{\dots}_{p\text{ entries}} ,\underbrace{0, 0, \dots,0,b}_{k+1\text{ entries}},\underbrace{\dots}_{q\text{ entries}});
\end{equation*}
we take zero the other entries of $G_{-(\ell+1)}$. Then the entries of $\wt B_{-(\ell)}$ outside of $(0, 0, \dots,0,b)$ are not changed.

  \item
Let $g_k$ be in the last column of $\un B_{1,\beta+1 }$  without its last entry. Then
\begin{equation*}\label{vbv1}
G_{-(\ell+1)}=(\,\underbrace{\text{any},
\emptyset}_{p\text{ entries}}\, ,\underbrace{g_1, g_2, \dots,g_k}_{k\text{ entries}}\,),\qquad\text{all }g_i\ne\emptyset,\ \ p\ge 0,
\end{equation*}
and \eqref{2340} takes the form
\begin{equation*}\label{fgh1}
\wt B_{-(\ell)}=B_{-(\ell)} +(\underbrace{\dots}_{p\text{ entries}} ,\underbrace{-g_1, g_1-g_2, \dots,g_{k-1}-g_k}_{k\text{ entries}}\,).
\end{equation*}
We choose $g_1, g_2, \dots,g_k$ such that
\begin{equation*}\label{vfb}
\wt B_{-(\ell)}=(\underbrace{\dots}_{p\text{ entries}} ,\underbrace{0, 0, \dots,0}_{k\text{ entries}}\,),
\end{equation*}
and take zero the other entries of $G_{-(\ell+1)}$.
\end{itemize}

Using these transformations, we consecutively reduce $B_{-(\ell)}$, starting at the top entry, to the form that is described in Theorem \ref{t1q0}(a). The diagonals under $B_{-(\ell)}$ in $\un B_{1,\beta+1}$ are not changed.

We apply this reduction to
all lower diagonals  of $\un B_{1,\beta+1}$ and to its main diagonal, and
obtain $\un B_{1,\beta+1}$ in which all entries on the main diagonal and under it are in the form described in Theorem \ref{t1q0}(a).

 \item[Step 2,] \emph{in which we reduce the entries of $\un B_{1,\beta+1}$ over its main diagonal.}
Let $\ell\ge 1$ be such that all diagonals of $\un B_{1,\beta+1}$ under its $\ell$th  \emph{upper diagonal}
\[
B_{(\ell)}:=(b_{1,\ell+1},
b_{2,\ell+2},\dots,b_{r_{\beta+1 }-\ell,r_{\beta+1}})
\]
have the form described in Theorem \ref{t1q0}(a). We use the transformation \eqref{1230} given by $\un G_{1,\beta+1}$ in which only the $(\ell-1)$st  upper diagonal $G_{(\ell-1)}=(g_1,g_2,
\dots,g_{r_{\beta+1}-\ell+1})$ is nonzero; its entries are arbitrary.  This transformation adds the vector
\[
(g_1-g_2,g_2-g_3,\dots,g_{r_{\beta +1}-\ell}-g_{r_{\beta+1}-\ell+1})
\]
to $B_{(\ell)}$; we
make $B_{(\ell)}=(0,\dots,0)$ preserving all the diagonals under it. We repeat this reduction until we obtain $\un B_{1,\beta+1}$ in which all entries over the main diagonal are zero. The obtained $\un B_{1,\beta+1}$ is in the form described in Theorem \ref{t1q0}(a),  which completes its proof.
\end{description}


\subsection{Proof of Theorem
\ref{t1q0}(b)}

Let $(W,B)$ and $(W,B')$ be matrix pairs of the form \eqref{b6s0} in which all blocks $B_i$ and $B_i'$ have the form \eqref{1a20}.

Suppose that $S^{-1}(W,B)S=(W,B')$ for some nonsingular $S$. We must prove that $B=B'$. The matrix $S$ has the form
\begin{equation*}\label{9iu0}
S=\mat{S_0&S_1&\ddots&S_{k-1}\\
&S_0&\ddots&\ddots\\
&&\ddots&S_1\\
0_n&&&S_0
},\qquad S_0\text{ is nonsingular.}
\end{equation*}
Since $S_0$ commutes with $F:=J_n(0)$, we have
\begin{equation}\label{web0}
 S_0=
\mat{a_0&a_1&\ddots&a_{n-1}\\
&a_0&\ddots&\ddots\\
&&\ddots&a_1\\0&&&a_0}
=a_0I_n+a_1F+a_2F^2 +\cdots+a_{n-1}F^{n-1}
\end{equation}
for some $a_0\ne 0,a_1,a_2,\ldots$ from the field $\ff$.

Reasoning by induction, we fix $v\in\{1,\dots,k-1\}$, assume that
\begin{equation}\label{csy0}
B_1=B_1',\ \dots,\ B_{v-1}=B_{v-1}'\qquad\text{if } v\ge 2,
\end{equation}
 and
prove that $B_v=B_v'$. For each $u=1,2,\dots$, we define the blocks $B_i^{(u)}$ and $B_i^{\prime(u)}$
by
\[
\mat{F^u&B_1^{(u)}&\ddots&B_v^{(u)}\\
&F^u&\ddots&\ddots\\
&&\ddots&B_{1}^{(u)}\\
0_n&&&F^u
}:=
\mat{F&B_{1}&\ddots&B_{v}\\
&F&\ddots&\ddots\\
&&\ddots&B_{1}\\
0_n&&&F
}^u
\]
and
\[
\mat{F^u&B_1^{\prime(u)}&
\ddots&B_v^{\prime(u)}\\
&F^u&\ddots&\ddots\\
&&\ddots&B_{1}^{\prime(u)}\\
0_n&&&F^u
}:=
\mat{F&B'_{1}&\ddots&B'_{v}\\
&F&\ddots&\ddots\\
&&\ddots&B'_{1}\\
0_n&&&F
}^u.
\]
The equalities \eqref{csy0} ensure that $B_1^{(u)}=B_1^{\prime(u)}$, \dots, $B_{v-1}^{(u)}=B_{v-1}^{\prime(u)}$.
Since $BS=SB'$, we have $B^uS=SB^{\prime u}$, and so
\[
\mat{F^u&\ddots
&B_{v}^{(u)}\\
&\ddots&\ddots\\
0&&F^u
}
\mat{S_0&\ddots&S_{v}\\
&\ddots&\ddots\\
0_n&&S_0
}
=
\mat{S_0&\ddots&S_{v}\\
&\ddots&\ddots\\
0_n&&S_0
}
\mat{F^u&\ddots
&B_{v}^{\prime (u)}\\
&\ddots&\ddots\\
0&&F^u
}.
\]
Equating the upper-right blocks,
we obtain
\begin{multline}\label{njhh0}
F^uS_v+B_{1}^{(u)}S_{v-1}+\dots
+B_{v-1}^{(u)}S_1+B_{v}^{(u)}S_0\\
=S_0B_{v}^{\prime (u)}+S_1B_{v-1}^{(u)}+\dots
+S_{v-1}B_{1}^{(u)}+S_vF^u.
\end{multline}

Since the traces of matrices satisfy  $\tr(XY)=\tr(YX)$ for all $X$ and $Y$, the equality \eqref{njhh0} implies that $\tr(B_{v}^{(u)}S_0)=
\tr(S_0B_{v}^{\prime (u)})$, hence
\begin{equation}\label{xla0}
\tr\bigl((B_{v}^{(u)}-B_{v}^{\prime (u)})S_0\bigr)=0.
\end{equation}

The blocks $B_{v}^{(u)}$ and $B_{v}^{\prime (u)}$ are represented in the form
\[
B_{v}^{(u)}=B_vF^{u-1}+FB_vF^{u-2}
+F^2B_vF^{u-3}+
\dots+F^{u-1}B_v+s(F,B_1,\dots,B_{v-1}),
\]
\[
B_{v}^{\prime (u)}=B'_vF^{u-1}+FB'_vF^{u-2}
+F^2B'_vF^{u-3}+
\dots+F^{u-1}B'_v+s(F,B_1,\dots,B_{v-1}),
\]
in which $s(F,B_1,\dots,B_{v-1})$ is a sum of products of matrices from the set $\{F,B_1,\dots,B_{v-1}\}$. Therefore,
\[
B_{v}^{(u)}-B_{v}^{\prime (u)}
=(B_v-B'_v)F^{u-1}
+F(B_v-B'_v)F^{u-2}
+\dots
+F^{u-1}(B_v-B'_v).
\]
By \eqref{xla0},
\[
\tr\bigl((B_v-B'_v)F^{u-1}S_0\bigr)
+\tr\bigl(F(B_v-B'_v)F^{u-2}S_0\bigr)
+\dots
+\tr\bigl(F^{u-1}(B_v-B'_v)S_0\bigr)=0.
\]
Since $\tr(XYZ)=\tr(YZX)$ for all square matrices of the same size,
\[
\tr\bigl((B_v-B'_v)F^{u-1}S_0\bigr)
+\tr\bigl((B_v-B'_v)F^{u-2}S_0F\bigr)
+\dots
+\tr\bigl((B_v-B'_v)S_0F^{u-1}\bigr)=0
\]
By \eqref{web0}, $F$ commutes with $S_0$. Hence,
$
u\tr\bigl((B_v-B'_v)F^{u-1}S_0\bigr)
=0.$ Since the characteristic of the field $\ff$ is zero, $u\ne 0$,
and so $\tr\bigl((B_v-B'_v)F^{u-1}S_0\bigr)=0$.

Substituting \eqref{web0}, we get
\begin{equation*}\label{vjh}
a_0\tr\bigl((B_v-B'_v)F^{u-1}\bigr)
+a_1\tr\bigl((B_v-B'_v)F^{u}\bigr)
+\dots+
a_{n-1}\tr\bigl((B_v-B'_v)
F^{u+n-2}\bigr)=0,
\end{equation*}
in which $a_0\ne 0$. We consecutively take  $u=n,n-1,\dots,1$ and obtain
\[
\tr\bigl((B_v-B'_v)F^{n-1}\bigr)=0,\
\tr\bigl((B_v-B'_v)F^{n-2}\bigr)=0,\
\dots,\
\tr\bigl((B_v-B'_v)F^0\bigr)=0.
\]
For each $n\times n$ matrix $X=[x_{ij}]$ and $0\le k<n$,
\[
\tr(XF^k)=x_{k+1,1}+x_{k+2,2}
+x_{k+3,3}+\cdots
\]
is the sum of entries of its $k$th lower diagonal.
Since $B_v$ and $B'_v$ are of the form \eqref{1a20}, we have $B_v=B'_v$.

\section*{Acknowledgements}
The authors would like to thank the referee who provided useful and detailed comments.
V. Futorny was supported by the CNPq (304467/2017-0) and the FAPESP (2018/23690-6).  V.V.~Sergeichuk is grateful to the University of S\~ ao Paulo for hospitality and to the FAPESP for financial support (2018/24089-4).

\end{document}